\documentclass[12pt,a4]{amsart}

\usepackage[all]{xy}
\usepackage{amscd}
\usepackage{amssymb}

\newtheorem{thm}{Theorem}[section]
\newtheorem{lem}[thm]{Lemma}
\newtheorem{cor}[thm]{Corollary}
\newtheorem{prop}[thm]{Proposition}

\theoremstyle{definition}

\newtheorem{rem}[thm]{Remark}
\newtheorem{defn}[thm]{Definition}

\newtheorem{ex}[thm]{Example}

\theoremstyle{remark}

\numberwithin{equation}{section}

\def\D{{\mathfrak D}}
\def\E{{\mathfrak E}}
\def\F{{\mathbb F}}

\def\Q{{\mathbb Q}}
\def\Z{{\mathbb Z}}
\def\C{{\mathbb C}}
\def\O{{\mathcal O}}

\def\T{{\mathfrak T}}
\def\U{{\mathfrak U}}
\def\p{{\mathfrak p}}
\def\q{{\mathfrak q}}

\def\X{{\mathfrak X}}
\def\Y{{\mathfrak Y}}
\def\ZZ{{\mathfrak Z}}

\def\Exc{\text{\rm Exc}}
\def\Frob{\text{\rm Frob}}
\def\Gr{\text{\rm Gr}}
\def\Gal{\text{\rm Gal}}
\def\Hom{\text{\rm Hom}}
\def\id{\text{\rm id}}
\def\GL{\text{\rm GL}}
\def\Spec{\text{\rm Spec}\,}
\def\Tr{\text{\rm Tr}}

\begin{document}

\title[Stringy Hodge numbers]
{Stringy Hodge numbers and p-adic Hodge theory}

\author[Tetsushi Ito]{Tetsushi Ito}
\address{Department of Mathematical Sciences,
University of Tokyo, 3-8-1 Komaba, Meguro, Tokyo 153-8914, Japan}
\email{itote2\char`\@ms.u-tokyo.ac.jp}

\address{Max-Planck-Institut F\"ur Mathematik,
Vivatsgasse 7, D-53111 Bonn, Germany}
\email{tetsushi\char`\@mpim-bonn.mpg.de}


\date{July 20, 2003}

\begin{abstract}
The aim of this paper is to give an application of
$p$-adic Hodge theory to stringy Hodge numbers
introduced by V. Batyrev for a mathematical formulation
of mirror symmetry.
Since the stringy Hodge numbers of an algebraic variety
are defined by choosing a resolution of singularities,
the well-definedness is not clear from the definition.
We give a proof of the well-definedness
by using arithmetic techniques such as $p$-adic integration
and $p$-adic Hodge theory.
Note that another proof of the well-definedness
was already obtained by V. Batyrev himself
by motivic integration.
\end{abstract}

\maketitle

\section{Introduction} 

Let $X$ be an irreducible normal algebraic variety over $\C$ with
at worst log-terminal singularities.
Let $\rho \colon Y \to X$
be a resolution of singularities
such that the exceptional divisor $\Exc(\rho)$ is
a normal crossing divisor
whose irreducible components $D_1,\ldots,D_r$ are smooth.
Let $K_Y = \rho^{\ast} K_X + \sum_{i=1}^{r} a_i D_i$
with $a_i \in \Q,\ a_i > -1$,
$I := \{ 1,\ldots,r \}$,
$D^{\circ}_J := \big( \bigcap_{j \in J} D_j \big)
  \backslash \big( \bigcup_{j \in I \backslash J} D_j \big)$
for a nonempty subset $J \subset I$, and
$D^{\circ}_{\emptyset} := Y \backslash \Exc(\rho)$.
We define the {\it stringy $E$-function} $E_{st}(X;u,v)$
of $X$ by the formula
$$ E_{st}(X;u,v) := \sum_{J \subset I} E(D^{\circ}_J;u,v)
     \prod_{j \in J} \frac{uv-1}{(uv)^{a_j+1}-1}, $$
where
$E(D^{\circ}_J;u,v) := \sum_{k} (-1)^k \sum_{i,j}
 h^{i,j}(\Gr_{i+j}^W H_c^k(D^{\circ}_J,\Q))\,u^i v^j$
is the generating function of the Hodge numbers of
$D^{\circ}_J$ (for details, see \S 2).
The aim of this paper is to give an alternative
proof of the following theorem by using
arithmetic techniques such as $p$-adic integration
and $p$-adic Hodge theory.

\begin{thm}[\cite{Batyrev2}, Theorem 3.4]
\label{MainTheorem}
The stringy $E$-function $E_{st}(X;u,v)$ defined as above
is independent of the choice of a resolution of singularities
$\rho \colon Y \to X$.
\end{thm}

Assume that $E_{st}(X;u,v)$ is a polynomial in $u,v$.
We define the {\it stringy Hodge numbers}
$h^{i,j}_{st}(X)$ of $X$ by the formula
$$ E_{st}(X;u,v) = \sum_{i,j} (-1)^{i+j}\,h^{i,j}_{st}(X)\,u^i v^j. $$
Therefore, by Theorem \ref{MainTheorem},
we establish the well-definedness of stringy Hodge numbers
(see also \cite{Batyrev2}).

Here we briefly recall a motivation of stringy Hodge numbers.
A mathematical formulation of mirror symmetry predicts
a symmetry between Hodge numbers of mirror varieties
(\cite{Batyrev2},\cite{Morrison}).
However, some examples discovered by physicists
show that the mirror of a smooth variety is not necessarily smooth.
In some cases, usual Hodge theory doesn't work well.
To overcome this difficulty, Batyrev introduced
stringy Hodge numbers as above (\cite{Batyrev2}).
Note that, for proper smooth varieties,
stringy Hodge numbers coincide with usual Hodge numbers
(Corollary \ref{Cor1}).
Today, several examples of stringy Hodge numbers are computed
from the viewpoint of mirror symmetry
(\cite{BatyrevBorisov},\cite{BatyrevDais},\cite{BorisovMavlyutov}).

There is an interesting history about
the proofs of Theorem \ref{MainTheorem}.
First of all, Batyrev proved that birational Calabi-Yau
manifolds have equal Betti numbers by using arithmetic techniques
such as $p$-adic integration and the Weil conjecture
(\cite{Batyrev1}).
Batyrev's method was generalized to birational
smooth minimal models by Wang (\cite{Wang1}).
The author obtained the equality of Hodge numbers
by using $p$-adic Hodge theory (\cite{Ito1},\cite{Ito2}).
Wang informed the author that he also obtained
the same result independently (\cite{Wang2},\cite{Wang3}).
On the other hand, in order to generalize Batyrev's work
on Betti numbers to Hodge numbers,
Kontsevich and Denef-Loeser developed
the theory of motivic integration,
which is a geometric analogue of $p$-adic integration
(\cite{Kontsevich},\cite{DenefLoeser2}).
Then, Batyrev introduced stringy Hodge numbers and
proved Theorem \ref{MainTheorem}
by using motivic integration (\cite{Batyrev2}).
In this paper, we give an alternative proof of
Theorem \ref{MainTheorem} by using arithmetic techniques
such as $p$-adic integration and $p$-adic Hodge theory.
In some sense, this paper goes back to
Batyrev's original arithmetic approach to
Theorem \ref{MainTheorem} by using $p$-adic Hodge theory.

It is worth mentioning that Theorem \ref{MainTheorem} has
non-trivial applications to birational geometry.
Firstly, for an algebraic variety $X$ over $\C$ with
a crepant resolution $\rho \colon Y \to X$,
the Hodge numbers of $Y$ are independent
of the choice of a crepant resolution $\rho \colon Y \to X$
(Corollary \ref{Cor2}).
This is an important fact in the study of McKay correspondences
in higher dimensions (\cite{BatyrevDais}).
Secondly, birational smooth minimal models
(e.g. Calabi-Yau manifolds) have equal Hodge numbers
(Corollary \ref{Cor3}).
Note that, in dimension $\leq 3$, this can also be proved by
the minimal model program (\cite{KMM},\cite{Kawamata},\cite{Kollar}).
Recently, a new proof valid in any dimension
was given by the weak factorization theorem
of birational maps (\cite{AKMW},\cite{Veys1},\cite{Veys2}).

This work is a continuation of author's previous works
(\cite{Ito1},\cite{Ito2}).
Here we note the new ingredients of this paper.
Basically, the main ideas are the same as before.
However, to treat stringy Hodge numbers rather than usual Hodge
numbers, we calculate some $p$-adic integration explicitly
(Proposition \ref{padicIntegration2}).
Furthermore, to treat combinations of cohomology groups of
open varieties, we generalize arithmetic techniques to open varieties
by a method of Deligne (\cite{HodgeI},\cite{HodgeII})
and work on the level of
a Grothendieck group of Galois representations
rather than individual cohomology groups (\S 5).

\vspace{0.2in}

\noindent
{\bf Acknowledgments.}
The author is grateful to Takeshi Saito and Kazuya Kato
for their advice and support.
He would like to thank Victor V. Batyrev,
Fran\c cois Loeser for helpful discussion on
motivic integration,
Shinichi Mochizuki, Takeshi Tsuji for
invaluable suggestion on $p$-adic Hodge theory,
Yujiro Kawamata,
Daisuke Matsushita, Yasunari Nagai,
Keiji Oguiso, Atsushi Takahashi, Hokuto Uehara,
Ezra Getzler, Willem Veys, Chin-Lung Wang
for comments and encouragement.
He would also like to thank Yoichi Mieda
for reading an earlier version carefully.
The author was supported by the Japan Society for the
Promotion of Science Research Fellowships for Young
Scientists.

\section{Stringy Hodge numbers}

In this section, we recall the definition of
stringy $E$-functions and stringy Hodge numbers as
in \cite{Batyrev2}.

For an algebraic variety $X$ over $\C$,
the cohomology groups with compact support $H_c^k(X,\Q)$
have canonical mixed Hodge structures by Deligne
(\cite{HodgeII}, \cite{HodgeIII}).
Let $W$ be the weight filtration on $H_c^k(X,\Q)$.
Each graded quotient $\Gr_l^W H_c^k(X,\Q)$ has a pure Hodge structure
of weight $l$. Let $h^{i,j}(\Gr_{i+j}^W H_c^k(X,\Q))$ be
the dimension of the $(i,j)$-th Hodge component of
$\Gr_{i+j}^W H_c^k(X,\Q)$ for each $i,j$.

\begin{defn}
\label{DefinitionE-function}
We define the {\it $E$-function} of $X$ as follows :
$$ E(X;u,v)
  := \sum_{k} (-1)^k \sum_{i,j} h^{i,j}(\Gr_{i+j}^W H_c^k(X,\Q))\,u^i v^j. $$
\end{defn}

For a proper smooth variety $X$ over $\C$,
$$ E(X;u,v) = \sum_{i,j} (-1)^{i+j} h^{i,j}(X) u^i v^j, $$
where $h^{i,j}(X) := \dim H^j(X,\Omega^i_X)$ are
the Hodge numbers of $X$ as usual.

\begin{rem}
$E$-function satisfies the following properties :
\begin{enumerate}
\item For $Z \subset X$, we have
  $ E(X;u,v) = E(X \backslash Z;u,v) + E(Z;u,v). $
\item For $X,Y$, we have
  $ E(X \times Y;u,v) = E(X;u,v) \cdot E(Y;u,v). $
\end{enumerate}
These two properties imply that it is natural to consider
$E$-function as a ring homomorphism from a Grothendieck group of
algebraic varieties over $\C$ to $\Z[u,v]$
(\cite{DenefLoeser1},\cite{DenefLoeser2},\cite{DenefLoeser3}).
\end{rem}

Let $X$ be an irreducible normal algebraic variety over $\C$,
and $\rho \colon Y \to X$ be a resolution of singularities such that
the exceptional divisor $\Exc(\rho)$ is a SNCD
(= simple normal crossing divisor).
Recall that a NCD (= normal crossing divisor) is
simple if its irreducible components are smooth.
Let the irreducible components of $\Exc(\rho)$
be $D_1,\ldots,D_r$.

\begin{defn}[\cite{Batyrev2}, Definition 2.2]
\label{DefLogTerminal}
$X$ is said to have {\it at worst log-terminal singularities}
if the following conditions are satisfied :
\begin{enumerate}
\item The canonical divisor $K_X$ is a $\Q$-Cartier divisor
(i.e. $X$ is $\Q$-Gorenstein).
\item We have
$$ K_Y = \rho^{\ast} K_X + \sum_{i=1}^{r} a_i D_i
   \qquad (a_i \in \Q), $$
with $a_i > -1$ (Note that this condition is independent of the choice
of a resolution $\rho \colon Y \to X$).
\end{enumerate}
\end{defn}

Let $X,Y$ be as above and $X$ have at worst log-terminal singularities.
Let $I := \{ 1,\ldots,r \}$.
For any subset $J \subset I$, we set
$$ D_J :=
\begin{cases} \bigcap_{j \in J} D_j & J \neq \emptyset \\
  Y & J = \emptyset
\end{cases},\qquad
D^{\circ}_J := D_J \backslash \bigcup_{j \in I \backslash J} D_j. $$

\begin{defn}[\cite{Batyrev2}, Definition 3.1]
\label{DefinitionStringyE-function}
We define the {\it stringy $E$-function} of $X$
as follows
$$ E_{st}(X;u,v) := \sum_{J \subset I} E(D^{\circ}_J;u,v)
     \prod_{j \in J} \frac{uv-1}{(uv)^{a_j+1}-1}, $$
where $E(D^{\circ}_J;u,v)$ is the $E$-function of a smooth variety
$D^{\circ}_J$ defined in the beginning of this section.
\end{defn}

\begin{rem}
Since $X$ has at worst log-terminal singularities,
$a_i + 1 > 0$ and hence the denominator of $E_{st}(X;u,v)$
doesn't vanish
(see also Remark \ref{p-adicIntegrationConvergence}).
In general, $E_{st}(X;u,v)$ is an element of
$\Q(u^{1/d},v^{1/d}) \cap \Z[[u^{1/d},v^{1/d}]]$,
where $d$ is the least common multiplier
of the denominators of $a_i$.
\end{rem}

\begin{defn}
Assume that $E_{st}(X;u,v)$ is a polynomial in $u,v$.
Then we define the {\it stringy Hodge numbers}
$h^{i,j}_{st}(X)$ of $X$ by the formula
$$ E_{st}(X;u,v) = \sum_{i,j} (-1)^{i+j}\,h^{i,j}_{st}(X)\,u^i v^j. $$
\end{defn}

Theorem \ref{MainTheorem} claims that
$E_{st}(X;u,v)$ is independent of the choice of a resolution
$\rho \colon Y \to X$.
Once we know the well-definedness,
we can prove some fundamental properties of $E_{st}(X;u,v)$
as in \cite{Batyrev2}.

Here we list some immediate corollaries of Theorem \ref{MainTheorem}.

\begin{cor}[\cite{Batyrev2}, Corollary 3.6]
\label{Cor1}
If $X$ is smooth, then we have $E_{st}(X;u,v) = E(X;u,v)$.
\end{cor}

\begin{proof}
This is clear because the identity map $\id \colon X \to X$
is a resolution of singularities.
\end{proof}

\begin{cor}[\cite{Batyrev2}, Theorem 3.12]
\label{Cor2}
Let $X$ be a projective algebraic variety over $\C$
which has a {\it crepant resolution} $\rho \colon Y \to X$
(i.e. $\rho \colon Y \to X$ is a resolution of singularities
with $\rho^{\ast} K_X = K_Y$).
Then the stringy Hodge numbers of $X$ are
equal to the Hodge numbers of $Y$ :
$$ h^{i,j}_{st}(X) = h^{i,j}(Y) \quad \text{for all}\ i,j. $$
In particular, the Hodge numbers of $Y$ are independent of
the choice of a crepant resolution $\rho \colon Y \to X$.
Moreover, we can compute the Hodge numbers of a crepant resolution
of $X$ via any resolution of singularities which is not necessarily
crepant.
\end{cor}

\begin{proof}
This is clear because we have
$E_{st}(X;u,v) = E(Y;u,v)$
by Definition \ref{DefinitionStringyE-function}.
\end{proof}

\begin{cor}
[\cite{Batyrev1},\cite{Batyrev2},\cite{Veys1},\cite{Veys2},\cite{Wang1},\cite{Wang2},\cite{Wang3},\cite{Ito1},\cite{Ito2}]
\label{Cor3}
Let $X,Y$ be projective smooth algebraic varieties over $\C$
whose canonical bundles are nef (i.e. $X,Y$ are {\em minimal models}).
Assume that $X,Y$ are birational.
Then $X,Y$ have equal Hodge numbers:
$$ h^{i,j}(X) = h^{i,j}(Y) \quad \text{for all}\ i,j. $$
\end{cor}

\begin{proof}
Let $f \colon X \dasharrow Y$ be a birational map.
Then we can find a projective smooth algebraic variety
$Z$ over $\C$ and birational morphisms
$g \colon Z \to X,\ h \colon Z \to Y$ such that
$f \circ g = h$ as birational maps
and $g^{\ast} K_X = h^{\ast} K_Y$.
$$
\xymatrix{
& Z \ar[dl]_{g} \ar[dr]^{h} \\
X \ar@{-->}[rr]^{f} & & Y }
$$
This is a standard fact in birational geometry
(for example, see \cite{Ito2}, Proposition 2.1).
We consider $g \colon Z \to X$ (resp. $h \colon Z \to Y$)
as a resolution of singularities of $X$ (resp. $Y$)
and calculate the stringy Hodge numbers of $X$ (resp. $Y$).
Since $g^{\ast} K_X = h^{\ast} K_Y$,
we have
$$ E(X;u,v) = E_{st}(X;u,v) = E_{st}(Y;u,v) = E(Y;u,v). $$
Hence we have the equality of the Hodge numbers of $X,Y$.
\end{proof}

\section{$p$-adic integration}

In this section, we recall Weil's $p$-adic integration
developed in \cite{Weil}
which is an important tool to count the number of
rational points valued in a finite field.

\subsection{Setup}

Let $p$ be a prime number and $\Q_p$ be the field of
$p$-adic numbers.
Let $F$ be a finite extension of $\Q_p$,
$R \subset F$ be the ring of integers in $F$,
$m \subset R$ be the maximal ideal of $R$,
$\F_q = R/m$ be the residue field of $F$ with $q$ elements,
where $q$ is a power of $p$.
For an element $x \in F$, we define
the {\it $p$-adic absolute value} $|x|_p$ by
$$ |x|_p :=
   \begin{cases} q^{-v(x)} & x \neq 0 \\ 0 & x=0 \end{cases} $$
where $v \colon F^{\times} \to \Z$ is
the normalized discrete valuation of $F$.

Let $\X$ be a smooth scheme over $R$ of
relative dimension $n$.
We can compute the number of the $\F_q$-rational points
$|\X(\F_q)|$ by integrating certain $p$-adic measure on
the set of $R$-rational points $\X(R)$.
We note that $\X(R)$ is a compact and totally disconnected
topological space with respect to its $p$-adic topology.

\subsection{$p$-adic integration of regular $n$-forms}

Let $\omega \in \Gamma(\X,\Omega^n_{\X/R})$ be
a {\it regular $n$-form} on $\X$, where $\Omega^n_{\X/R}$ is
the relative canonical bundle of $\X/R$.
We shall define the $p$-adic integration of $\omega$ on $\X(R)$
as follows.
Let $s \in \X(R)$ be a $R$-rational point.
Let $U \subset \X(R)$ be a sufficiently small
$p$-adic open neighborhood of $s$
on which there exists a system of local $p$-adic
coordinates $\{ x_1,\ldots,x_n \}$.
Then $\{ x_1,\ldots,x_n \}$ defines a $p$-adic analytic map
$$ x = (x_1,\ldots,x_n) \colon U \longrightarrow R^n $$
which is a homeomorphism between $U$ and
a $p$-adic open set $V$ of $R^n$.
By using the above coordinates, $\omega$ is written as
$$ \omega = f(x) \ dx_1 \wedge \cdots \wedge dx_n. $$
We consider $f(x)$ as a $p$-adic analytic function on $V$.
Then we define the $p$-adic integration of $\omega$ on $U$
by the equation
$$ \int_{U} |\omega|_p := \int_{V} |f(x)|_p \ dx_1 \cdots dx_n, $$
where $|f(x)|_p$ is the $p$-adic absolute value of the value of $f$
at $x \in V$ and $dx_1 \cdots dx_n$ is the Haar measure on $R^n$
normalized by the condition
$$ \int_{R^n} dx_1 \cdots dx_n = 1. $$
By patching them, we get
the {\it $p$-adic integration of $\omega$ on $\X(R)$}
$$ \int_{\X(R)} |\omega|_p. $$

\subsection{$p$-adic integration of gauge forms}

By definition, a {\it gauge form} $\omega$ on $\X$ is
a nowhere vanishing global section
$\omega \in \Gamma(\X,\Omega^n_{\X/R})$.
The most important property of $p$-adic integration is that
the $p$-adic integration of a gauge form
computes the number of $\F_q$-rational points.

\begin{prop}[\cite{Weil}, 2.2.5]
\label{padicIntegration1}
Let $\X$ be a smooth scheme over $R$ of relative dimension $n$
and $\omega$ be a gauge form on $\X$.
Then
$$ \int_{\X(R)} |\omega|_p = \frac{|\X(\F_q)|}{q^n}. $$
\end{prop}

\begin{proof}
Let
$$ \varphi \colon \X(R) \longrightarrow \X(\F_q) $$
be the reduction map.
For $\bar{x} \in \X(\F_q)$, $\varphi^{-1}(\bar{x})$ is
a $p$-adic open set of $\X(R)$.
Therefore, it is enough to show
$$ \int_{\varphi^{-1}(\bar{x})} |\omega|_p = \frac{1}{q^n}. $$
Let $\{ x_1,\ldots,x_n \} \subset \O_{\X,\bar{x}}$ be
a regular system of parameters at $\bar{x}$.
Then $\{ x_1,\ldots,x_n \}$ defines a system of local
$p$-adic coordinates on $\varphi^{-1}(\bar{x})$ and
$$ x = (x_1,\ldots,x_n) \colon \varphi^{-1}(\bar{x})
      \longrightarrow m^n \subset R^n $$
is a $p$-adic analytic homeomorphism.
Let $\omega$ be written as
$\omega = f(x) \ dx_1 \wedge \cdots \wedge dx_n$.
Since $\omega$ is a gauge form, $f(x)$ is a $p$-adic unit
for all $x \in \varphi^{-1}(\bar{x})$.
Therefore $|f(x)|_p = 1$.
Then we have
$$ \int_{\varphi^{-1}(\bar{x})} |\omega|_p
     = \int_{m^n} dx_1 \cdots dx_n = \frac{1}{q^n} $$
since $m^n$ is an index $q^n$ subgroup of $R^n$.
\end{proof}

\subsection{Computation of some $p$-adic integration}

Since we treat stringy Hodge numbers rather than usual
Hodge numbers in this paper,
we need to compute some $p$-adic integration slightly
more general than Proposition \ref{padicIntegration1}.

Firstly, we generalize $p$-adic integration
to an $r$-pluricanonical form with pole
for $r \in \Z,\ r \geq 1$ (\cite{Wang1}).
An {\it $r$-pluricanonical form} on $\X$ is a section of
$(\Omega^n_{\X/R})^{\otimes r}$ over $\X$.
An {\it $r$-pluricanonical form with pole} on $\X$
is a section of $(\Omega^n_{\U/R})^{\otimes r}$ over $\U$
for some open subscheme $\U \subset \X$.
Let $\omega$ be an $r$-pluricanonical form with pole on $\X$.
As in the case of a regular $n$-form,
locally in $p$-adic topology, $\omega$ is written as
$$ \omega = f(x) \ (dx_1 \wedge \cdots \wedge dx_n)^{\otimes r} $$
for a system of local $p$-adic coordinates $\{ x_1,\dots,x_n \}$.
Note that $f(x)$ is a $p$-adic analytic function with pole.
Then we put
$$ \int_{U} |\omega|_p^{1/r} := 
   \int_{V} |f(x)|_p^{1/r} \ dx_1 \cdots dx_n, $$
where $U,V$ are the same as in the case of a regular $n$-form,
if the right hand side converges.
If the above integral converges for each open neighborhood,
by patching them, 
we get the {\it $p$-adic integration of an $r$-pluricanonical form
with pole $\omega$ on $\X(R)$}
$$ \int_{\X(R)} |\omega|_p^{1/r}. $$
Note that, if $\omega$ has no pole,
the above integral always converges.

\begin{rem}
For an $r$-pluricanonical form with pole $\omega$,
$\omega^{\otimes s}$ is an $rs$-pluricanonical form with pole.
If the $p$-adic integration of $\omega$ converges,
then the $p$-adic integration of $\omega^{\otimes s}$ also
converges and they are equal :
$$ \int_{\X(R)} |\omega|_p^{1/r}
     = \int_{\X(R)} |\omega^{\otimes s}|_p^{1/rs}, $$
\end{rem}

Before computing $p$-adic integration, we recall
the notion of relative SNCD.

\begin{defn}
\label{DefinitionRelativeSNCD}
Let $f \colon \X \to S$ be a proper smooth morphism of schemes
and $\D = \sum_{i=1}^{r} a_i \D_i\ (a_i \in \Q,\ a_i \neq 0)$
be a $\Q$-divisor on $\X$.
Let $\text{Supp}\,\D = \bigcup_{i=1}^{r} \D_i$ be
the support of $\D$.
Then, $\D$ is called a {\it relative SNCD}
(=simple normal crossing divisor) on $\X/S$
if all $\D_i$ are smooth over $S$ and,
for all $x \in \text{Supp}\,\D$,
the completion of $\text{Supp}\,\D \hookrightarrow \X$ at $x$
is isomorphic to
$$
\Spec(\widehat{\O}_{S,f(x)}[[x_1,\ldots,x_d]]/(x_1 \cdots x_s))
\hookrightarrow
\Spec(\widehat{\O}_{S,f(x)}[[x_1,\ldots,x_d]])
$$
for some $s\ (1 \leq s \leq d)$,
where $d$ is the relative dimension of $f$.
\end{defn}

Note that
$\widehat{\O}_{\X,x}$ is isomorphic to
$\widehat{\O}_{S,f(x)}[[x_1,\ldots,x_d]]$
because $f$ is smooth of relative dimension $d$.
In this case, for a nonempty subset $J \subset \{ 1,\ldots,r \}$,
$\bigcap_{j \in J} \D_j$ is smooth of relative dimension
$d - |J|$ over $S$.

We shall compute some $p$-adic integration.
Let $\X$ be a smooth scheme over $R$ of relative dimension $n$,
and $\omega$ an $r$-pluricanonical form with pole on $\X$.
Assume that
$$ \text{\rm div}(\omega) = \sum_{i=1}^{s} a_i \D_i $$
is a relative SNCD on $\X/R$.
Let $I := \{ 1,\ldots,s \}$.
For any subset $J \subset I$, we set
$$ \D_J :=
\begin{cases} \bigcap_{j \in J} \D_j & J \neq \emptyset \\
  \X & J = \emptyset
\end{cases},\qquad
\D^{\circ}_J := \D_J \backslash \bigcup_{j \in I \backslash J} \D_j. $$

\begin{prop}
\label{padicIntegration2}
Let notation be as above.
If $a_i > -r$ for all $i \in I$,
then the $p$-adic integration of $\omega$ on $\X(R)$
converges, and we have the following equality :
$$ \int_{\X(R)} |\omega|_p^{1/r}
     = \frac{1}{q^n} \sum_{J \subset I} |\D^{\circ}_J(\F_q)|
         \prod_{j \in J} \frac{q-1}{q^{(a_j/r)+1}-1}. $$
\end{prop}

If $r=1$ and $\omega$ is a gauge form,
Proposition \ref{padicIntegration2} is nothing but
Proposition \ref{padicIntegration1}.

\begin{proof}
The idea of proof is the same as in Proposition \ref{padicIntegration1}.
Let
$$ \varphi \colon \X(R) \longrightarrow \X(\F_q) $$
be the reduction map.
For $\bar{x} \in \X(\F_q)$, $\varphi^{-1}(\bar{x})$ is
a $p$-adic open set of $\X(R)$.
Therefore, it is enough to show
$$ \int_{\varphi^{-1}(\bar{x})} |\omega|_p^{1/r}
   = \frac{1}{q^n} \prod_{j \in I\ \text{s.t.}\ \bar{x} \in \D_j(\F_q)}
       \frac{q-1}{q^{(a_j/r)+1}-1}. $$
Let $\{ j_1,\ldots,j_k \} = \{ j \in I \mid \bar{x} \in \D_j(\F_q) \}$.
Let $\{ x_1,\ldots,x_n \} \subset \O_{\X,\bar{x}}$ be
a regular system of parameters at $\bar{x}$ such that
$\D_{j_i}$ is defined by $x_{i} = 0$ at $\bar{x}$ for all $i=1,\ldots,k$.
Then $\{ x_1,\ldots,x_n \}$ defines a system of local
$p$-adic coordinates on $\varphi^{-1}(\bar{x})$ and
$$ x = (x_1,\ldots,x_n) \colon \varphi^{-1}(\bar{x})
      \longrightarrow m^n \subset R^n $$
is a $p$-adic analytic homeomorphism.
Here $\omega$ is written as
$$ \omega = f(x) \cdot {x_1}^{\!a_{j_1}} \cdots {x_k}^{\!\!a_{j_k}}
              (dx_1 \wedge \cdots \wedge dx_n)^{\otimes r}, $$
where $f(x)$ is a $p$-adic unit for all $x \in \varphi^{-1}(\bar{x})$.
Hence we have
$$ |f(x) \cdot {x_1}^{\!a_{j_1}} \cdots {x_k}^{\!\!a_{j_k}}|_p^{1/r}
      = |x_1|_p^{a_{j_1}/r} \cdots\ |x_k|_p^{a_{j_k}/r}, $$
and
\begin{align*}
  \int_{\varphi^{-1}(\bar{x})} |\omega|_p^{1/r}
    &= \int_{m^n} |x_1|_p^{a_{j_1}/r} \cdots\ |x_k|_p^{a_{j_k}/r}
         dx_1 \cdots dx_n \\
    &= \frac{1}{q^{n-k}} \left(
         \int_{m^k} |x_1|_p^{a_{j_1}/r} \cdots\ |x_k|_p^{a_{j_k}/r}
         dx_1 \cdots dx_k \right).
\end{align*}
Therefore, it is enough to prove the following lemma.
\end{proof}

\begin{lem}
\label{padicIntegrationLemma}
For $k_1,\ldots,k_n \in \Q,\ k_i > -1$,
the $p$-adic integration
$$ \int_{m^n} |x_1|_p^{k_1} \cdots |x_n|_p^{k_n} dx_1 \cdots dx_n $$
converges and is equal to
$$ \frac{1}{q^n} \prod_{i=1}^{n} \frac{q-1}{q^{k_i+1}-1}. $$
\end{lem}

\begin{proof}
By iterated integration, we have
$$ \int_{m^n} |x_1|_p^{k_1} \cdots |x_n|_p^{k_n} dx_1 \cdots dx_n
    = \left( \int_{m} |x_1|_p^{k_1} dx_1 \right)
        \cdots \left( \int_{m} |x_n|_p^{k_n} dx_n \right). $$
Therefore, it is enough to prove
$$ \int_{m} |x|_p^{k}\,dx
     = \frac{1}{q} \cdot \frac{q-1}{q^{k+1}-1} $$
for $k \in \Q,\ k > -1$.

We compute the above integration by dividing $m$ as a disjoint union
of open subsets as follows
$$ m = \coprod_{i=1}^{\infty} m^i \backslash m^{i+1}. $$
For $x \in m^i \backslash m^{i+1}$, $|x|_p^k = q^{-ki}$.
The volume of $m^i$ is $q^{-i}$ with respect to
the normalized Haar measure on $R$ since $m^i$ is
an index $q^i$ subgroup of $R$.
Therefore, we have
\begin{align*}
 \int_{m} |x|_p^{k}\,dx
  &= \sum_{i=1}^{\infty} q^{-ik}\,\text{vol}\,(m^i \backslash m^{i+1})
   = \sum_{i=1}^{\infty} q^{-ik} ( q^{-i} - q^{-(i+1)} ) \\
  &= (1 - q^{-1}) \sum_{i=1}^{\infty} (q^{-(k+1)})^i.
\end{align*}
Since $k > -1$, this infinite sum converges to
$$ (1 - q^{-1}) \cdot \frac{q^{-(k+1)}}{1-q^{-(k+1)}}
     = \frac{1}{q} \cdot \frac{q-1}{q^{k+1}-1}. $$
Hence we have Lemma \ref{padicIntegrationLemma}
and the proof of Proposition \ref{padicIntegration2} is completed.
\end{proof}

\begin{rem}
A curious reader may notice the similarity
between the expression in Proposition \ref{padicIntegration2}
and Definition \ref{DefinitionStringyE-function}.
This is the starting point of our proof of
Theorem \ref{MainTheorem}.
However, to recover information of the Hodge numbers
of an algebraic variety from the numbers of rational points
valued in finite fields,
we need some deep arithmetic results as in \S 4, \S 5.
\end{rem}

\begin{rem}
\label{p-adicIntegrationConvergence}
As we easily see in the proof of Lemma \ref{padicIntegrationLemma},
the $p$-adic integration
$$ \int_{m^n} |x_1|_p^{k_1} \cdots |x_n|_p^{k_n} dx_1 \cdots dx_n $$
doesn't converge if $k_i \leq -1$ for some $i$.
This is the reason why we assume singularities are
at worst log-terminal in Theorem \ref{MainTheorem}.
\end{rem}

\section{Local Galois representations}

In this section, we recall some results on Galois representations
over a $p$-adic field.

\subsection{Setup}

Let $K$ be a number field.
Let $\p$ be a maximal ideal of $\O_K$.
Let $K_{\p}$ be a $\p$-adic completion of $K$,
$\O_{K_{\p}}$ the ring of integers of $K_{\p}$,
$\F_q = \O_{K}/\p$ the residue field of $K_{\p}$
with $q$ elements, and
$\overline{K}_{\p}$ (resp. $\overline{\F}_q$)
an algebraic closure of $K_{\p}$ (resp. $\F_q$).

We have an exact sequence
$$
\begin{CD}
0 @>>> I_{K_{\p}} @>>> \Gal(\overline{K}_{\p}/K_{\p})
  @>>> \Gal(\overline{\F}_q/\F_q) @>>> 0,
\end{CD}
$$
where $I_{K_{\p}}$ is called the {\it inertia group} at $\p$.
$\Gal(\overline{\F}_q/\F_q)$ is topologically generated by
the $q$-th power Frobenius automorphism $x \mapsto x^q$
of $\overline{\F}_q$.
The inverse of this automorphism is called
the {\it geometric Frobenius element} at $\p$
and denoted by $\Frob_{\p}$.

Let $X$ be a proper smooth variety over $K_{\p}$,
$l$ be a prime number, and $k$ be an integer.
Then the absolute Galois group $\Gal(\overline{K}_{\p}/K_{\p})$
acts continuously on the $l$-adic \'etale cohomology group
$H^k_{\text{\'et}}(X_{\overline{K}_{\p}}, \Q_l)$ of
$X_{\overline{K}_{\p}} = X \otimes_{K_{\p}} \overline{K}_{\p}$.
In the followings,
we recall some results on this
$\Gal(\overline{K}_{\p}/K_{\p})$-representation in two cases.

\subsection{The Weil conjecture}

Firstly, we assume that $\p$ doesn't divide $l$
and there exists a proper smooth scheme $\X$ over $\O_{K_{\p}}$
such that $\X \otimes_{\O_{K_{\p}}} K_{\p} = X$
($\X$ is called a {\it proper smooth model} of $X$
over $\O_{K_{\p}}$).

In this case, the action of $I_{K_{\p}}$ on
$H^k_{\text{\'et}}(X_{\overline{K}_{\p}}, \Q_l)$ is trivial
(i.e. the action of $\Gal(\overline{K}_{\p}/K_{\p})$ is
{\it unramified})
by the proper smooth base change theorem on \'etale cohomology.
Therefore, the action of $\Gal(\overline{K}_{\p}/K_{\p})$
is determined by the action of $\Frob_{\p}$.

By the Lefschetz trace formula for \'etale cohomology,
we have
$$ |\X(\F_q)| = \sum_{k} (-1)^k
   \,\Tr(\Frob_{\p};H^k_{\text{\'et}}(X_{\overline{K}_{\p}},\Q_l)). $$
Furthermore, the characteristic polynomial
$ P_k(t) = \det(1 - t \cdot \Frob_{\p};
    H^k_{\text{\rm \'et}}(X_{\overline{K}_{\p}}, \Q_l)) $
has integer coefficients and
all complex absolute values of all conjugates of
the roots of $P_k(t)$ are equal to $q^{-k/2}$ by
the Weil conjecture proved by Deligne (\cite{WeilI},\cite{WeilII}).

\subsection{$p$-adic Hodge theory}

Secondly, we assume $\p$ divides $l$.
Let $p = l$ in this subsection to avoid confusion.
Here no assumption is required for
a model of $X$ over $\O_{K_{\p}}$
(see Remark \ref{Remarkp-adicHodgeTheory}).
In this case,
the action of the inertia group $I_{K_{\p}}$ on
$H^k_{\text{\'et}}(X_{\overline{K}}, \Q_p)$ is
highly nontrivial.

Let $\C_p$ be a $\p$-adic completion of $\overline{K}_{\p}$
on which $\Gal(\overline{K}_{\p}/K_{\p})$ acts
continuously.
We recall the {\it Tate twists}.
Let
$\Q_p(0) := \Q_p,
\ \Q_p(1) := \left( \varprojlim \mu_{p^n} \right) \otimes_{\Z_p} \Q_p$.
For $n \geq 1$, let
$\Q_p(n) := \Q_p(1)^{\otimes n},\ \Q_p(-n) := \Hom(\Q_p(n),\Q_p)$.
Moreover, for a $\Gal(\overline{K}_{\p}/K_{\p})$-representation $V$
over $\Q_p$, we define $V(n) := V \otimes_{\Q_p} \Q_p(n)$
on which $\Gal(\overline{K}_{\p}/K_{\p})$ acts diagonally.

$p$-adic Hodge theory claims the {\it Hodge-Tate decomposition}
of $X$ as follows :
$$
 \bigoplus_{i,j\ \text{s.t.}\ i+j=k}
    H^j(X,\Omega_{X}^i) \otimes_{K} \C_p(-i)
 \cong
    H^k_{\text{\rm \'et}}(X_{\overline{K}_{\p}},\Q_p)
    \otimes_{\Q_p} \C_p.
$$
This is an isomorphism of
$\Gal(\overline{K}_{\p}/K_{\p})$-representations,
where $\Gal(\overline{K}_{\p}/K_{\p})$ acts on $H^i(X,\Omega_{X}^j)$
trivially and on the right hand side diagonally.
This is a $p$-adic analogue of the usual Hodge decomposition
over $\C$.

As a consequence, we can recover
the Hodge numbers of $X$ from its $p$-adic Galois representations
as follows :
$$ \dim_{K_{\p}} H^j(X,\Omega_{X}^i)
     = \dim_{K_{\p}}
       (H^{i+j}_{\text{\rm \'et}}(X_{\overline{K}_{\p}},\Q_p)
         \otimes_{\Q_p} \C_p(i))^{\Gal(\overline{K}_{\p}/K_{\p})}, $$
since $(\C_p)^{\Gal(\overline{K}_{\p}/K_{\p})} = K_{\p}$ and
$(\C_p(i))^{\Gal(\overline{K}_{\p}/K_{\p})} = 0$ for all $i \neq 0$
(\cite{Tate}, Theorem 2).

\begin{rem}
\label{Remarkp-adicHodgeTheory}
A proof of Hodge-Tate decomposition was given
by Faltings (\cite{Faltings} for recent developments
of Faltings' theory of almost \'etale extensions,
see also \cite{Faltings2}).
Tsuji gave another proof by reducing to the semi-stable reduction case
by de Jong's alteration (\cite{Tsuji}).
However, in this paper, we don't need the full version of
the Hodge-Tate decomposition.
For example, 
the result of Fontaine-Messing is enough for us (\cite{FontaineMessing}).
They proved the Hodge-Tate decomposition
when $K_{\p}$ is unramified over $\Q_p$,
$\dim X < p$, and $X$ has a proper smooth model over $\O_{K_{\p}}$.
\end{rem}

\begin{rem}
\label{RemarkSemisimplification}
Moreover, by Lemma \ref{LemmaSemisimplification} below,
we see that
the semisimplification of the $p$-adic Galois representation
determines the Hodge numbers by the same formula :
$$ \dim_{K_{\p}} H^j(X,\Omega_{X}^i)
     = \dim_{K_{\p}}
       (H^{i+j}_{\text{\rm \'et}}(X_{\overline{K}_{\p}},\Q_p)^{ss}
         \otimes_{\Q_p} \C_p(i))^{\Gal(\overline{K}_{\p}/K_{\p})}, $$
where {\it ss} denotes the semisimplification as
a $\Gal(\overline{K}_{\p}/K_{\p})$-representation.
This is a simple but important observation to consider
the Hodge-Tate decomposition on the level of
a Grothendieck group of Galois representations in \S 5. 
\end{rem}

\begin{lem}[\cite{Ito2}, Lemma 4.4]
\label{LemmaSemisimplification}
Let $0 \to V_1 \to V_2 \to V_3 \to 0$
be an exact sequence of finite dimensional
$\Gal(\overline{K}_{\p}/K_{\p})$-representations over $\Q_p$.
We define
$h^{n}(V_i) := \dim (V_i \otimes_{\Q_p} \C_p(n))^{
\Gal(\overline{K}_{\p}/K_{\p})}$
for $i = 1,2,3$ and an integer $n$.
Assume that $\dim V_2 = \sum_n h^n(V_2)$.
Then, we have
$\dim V_1 = \sum_n h^n(V_1)$,\ $\dim V_3 = \sum_n h^n(V_3)$
and $h^n(V_2) = h^n(V_1) + h^n(V_3)$
for all $n$.
\end{lem}

\begin{proof}
This lemma seems well-known to specialists.
However, we write the proof for reader's convenience.
In general, we have an inequality $\sum_n h^n(V_i) \leq \dim V_i$
for $i = 1,3$ (for example, see \cite{Fontaine}).
We shall prove these inequalities are in fact equalities.
Since the functor taking
$\Gal(\overline{K}_{\p}/K_{\p})$-invariant is left exact,
$$
\begin{CD}
0 @>>> (V_1 \otimes_{\Q_p} \C_p(n))^{\Gal(\overline{K}_{\p}/K_{\p})}
@>>> (V_2 \otimes_{\Q_p} \C_p(n))^{\Gal(\overline{K}_{\p}/K_{\p})} \\
@>>> (V_3 \otimes_{\Q_p} \C_p(n))^{\Gal(\overline{K}_{\p}/K_{\p})}
\end{CD}
$$
is exact. Therefore
$h^n(V_2) \leq h^n(V_1) + h^n(V_3)$ for all $n$.
Then we have
\begin{align*}
\dim V_2 &= \sum_n h^n(V_2)
         \leq \sum_n h^n(V_1) + \sum_n h^n(V_3)
         \leq \dim V_1 + \dim V_3 \\
         &= \dim V_2
\end{align*}
and hence Lemma \ref{LemmaSemisimplification}.
\end{proof}

\section{Global Galois representations}

In this section, we recall some results on Galois representations
over a number field.

\subsection{An application of the Chebotarev density theorem}

The following proposition is very important to work on
the level of a Grothendieck group of Galois representations
over a number field.
This is an application of the Chebotarev density theorem
in algebraic number theory.

\begin{prop}[\cite{Serre}, I.2.3]
\label{CorChebotarev}
Let $K$ be a number field and $l$ be a prime number.
Let $V,V'$
be two continuous $l$-adic $\Gal(\overline{K}/K)$-representations
such that they are unramified outside a finite
set $S$ of maximal ideals of $\O_K$
and satisfy
$$ \Tr(\Frob_{\p};V) = \Tr(\Frob_{\p};V')
   \qquad \text{for all} \quad \p \notin S. $$
Then $V$ and $V'$ have the same semisimplifications as
$\Gal(\overline{K}/K)$-representations.
\end{prop}

\begin{proof}
We only sketch the proof (for details, see \cite{Serre}).
By the representation theory of a group over
a field of characteristic 0
(see, for example, Bourbaki, {\it Alg\`ebre},
Ch. 8, \S 12, n${}^{\circ}$ 1, Prop 3.),
the semisimplification of
a $\Gal(\overline{K}/K)$-representation is determined
by the traces of all elements in $\Gal(\overline{K}/K)$.
Roughly speaking, the Chebotarev density theorem claims
the set of conjugates of $\Frob_{\p}$ for $\p \notin S$
is dense in $\Gal(\overline{K}/K)$.
Since $V$ and $V'$ are continuous representations,
the equality of the traces of all $\Frob_{\p}$ for $\p \notin S$
implies the equality of traces of all elements
in $\Gal(\overline{K}/K)$.
Hence we have Proposition \ref{CorChebotarev}.
\end{proof}

\subsection{Some Grothendieck groups of Galois representations}

Let $K$ be a number field.
Let $S$ be a finite set of maximal ideals of $\O_K$.
We fix a prime number $l = p$ and a maximal ideal $\p$
dividing $p$.
For every maximal ideal $\q$ of $\O_K$,
we fix an inclusion $\overline{K} \hookrightarrow \overline{K}_{\q}$.
Then we consider
$\Gal(\overline{K}_{\q}/K_{\q}) \subset \Gal(\overline{K}/K)$
for all $\q$.

\begin{defn}
\label{DefinitionGrothendieckGroups}
Let $K(l,S,\p)$ be an abelian group
generated by $\Gal(\overline{K}/K)$-representations $V$
satisfying the following conditions
modulo an equivalence relation $\sim$ generated by
$[V_1] + [V_3] \sim [V_2]$
for an exact sequence $0 \to V_1 \to V_2 \to V_3 \to 0$ :
\begin{enumerate}
\item ({\it unramifiedness} outside $S$)
\ For $\q \notin S$, $I_{K_{\q}}$ acts on $V$ trivially.
\item ({\it weight filtration} outside $S$)
\ There exists a unique increasing 
$\Gal(\overline{K}/K)$-equivariant filtration $W$ on $V$
indexed by integers satisfying the following
conditions :
\begin{enumerate}
\item $W_k V = 0$ for $k \ll 0$, $W_k V = V$ for $k \gg 0$.
\item For every integer $k$ and $\q \notin S$,
the characteristic polynomial
$P_k(t) = \det(1 - t \cdot \Frob_{\q};\Gr_k^W V)$
has integer coefficients and
all complex absolute values of all conjugates of
the roots of $P_k(t)$ are equal to $|\O_K/{\q}|^{-k/2}$.
\end{enumerate}
$W$ is called the {\it weight filtration} on $V$.
\item ({\it Hodge-Tate decomposition} at $\p$)
\ For integers $i,j$, we define
$$
h_{\p}^{i,\,j}(V) := \dim_{K_{\p}}
(\Gr_{i+j}^W V \otimes_{\Q_p} \C_p(i))^{\Gal(\overline{K}_{\p}/K_{\p})}.
$$
Then these numbers satisfy
$ \sum_{i, j} h_{\p}^{i,j}(V) = \dim_{\Q_p} V. $
\end{enumerate}
$K(l,S,\p)$ is called the {\it Grothendieck group} of
$l$-adic $\Gal(\overline{K}/K)$-representations
which are unramified outside $S$, have weight filtration,
and have Hodge-Tate decomposition at $\p$.
Let $[V]$ denote the class of $V$ in $K(l,S,\p)$.
An element in $K(l,S,\p)$ is called a
{\it virtual $\Gal(\overline{K}/K)$-representation}.
\end{defn}

Similarly, for an integer $k$, we define
$K(l,S,\p,k)$ as a subgroup of $K(l,S,\p)$ generated by
$V$ satisfying $\Gr_k^W V = V$.
$K(l,S,\p,k)$ is called the Grothendieck group of
$l$-adic $\Gal(\overline{K}/K)$-representations
which are unramified outside $S$, have {\it weight $k$},
and have Hodge-Tate decomposition at $\p$.

By the Jordan-H\"older theorem,
$K(l,S,\p)$ is a free abelian group generated by
simple $\Gal(\overline{K}/K)$-representations
in $K(l,S,\p)$.
Since $[V] = \sum_k [\Gr_k^W V]$ in $K(l,S,\p)$,
a simple $\Gal(\overline{K}/K)$-representation has
only one weight.
Therefore we have a direct sum decomposition as follows :
$$ K(l,S,\p) = \bigoplus_{k \in \Z} K(l,S,\p,k). $$

We define a ring structure
on $K(l,S,\p)$ by extending the tensor product
$[V_1] \cdot [V_2] = [V_1 \otimes V_2]$.
Then $K(l,S,\p)$ has a structure of a graded ring
by the direct sum decomposition as above.

\begin{defn}
\label{Definitionp-adicE-function}
For a $\Gal(\overline{K}/K)$-representation $V$ in $K(l,S,\p)$,
we define the {\it $p$-adic $E$-function} of $V$ as follows
$$ E_{\p}(V;u,v)
     := \sum_{i,j} h_{\p}^{i,j}(V)\,u^i v^j. $$
\end{defn}

\begin{rem}
It is easy to see that $p$-adic $E$-function satisfies
the following properties
(see Remark \ref{RemarkSemisimplification}) :
\begin{enumerate}
\item For an exact sequence $0 \to V_1 \to V_2 \to V_3 \to 0$,
we have $E_{\p}(V_1;u,v) + E_{\p}(V_3;u,v) = E_{\p}(V_2;u,v)$.
\item For $V_1,V_2$, we have
$E_{\p}(V_1 \otimes V_2;u,v) = E_{\p}(V_1;u,v) \cdot E_{\p}(V_2;u,v)$.
\end{enumerate}
Therefore, we can extend $p$-adic $E$-function
to a ring homomorphism
$$ E_{\p} \colon K(l,S,\p) \longrightarrow \Z[u,v]. $$
\end{rem}

\subsection{A variant --- Galois represenations
with fractional weight filtration}

Let $d$ be an integer.
Here we introduce a variant of $K(l,S,\p)$ whose
weight filtration is indexed by elements of $\frac{1}{d} \Z$
instead of $\Z$.
This generalization is necessary to treat
the $\Q$-Gorenstein case in the proof of
Theorem \ref{MainTheorem}.

Firstly, we introduce the {\it fractional Tate twists}
$\Q_p(a) \ (a \in \frac{1}{d}\Z)$ as follows
(for usual Tate twists, see \S 4.3).
Let $L$ be a field and $p$ be a prime number.
$\Q_p(1)$ is a one dimensional
$\Gal(\overline{L}/L)$-representation
$$ \rho \colon \Gal(\overline{L}/L) \longrightarrow
     \GL(\Q_p(1))
     \cong \GL(1,\Q_p) = \Q_p^{\times} $$
whose image is contained in $\Z_p^{\times}$.
There exist open subgroups $U \subset \Z_p^{\times}$,
$V \subset \Z_p$ on which
$\log \colon U \to V$ and $\exp \colon V \to U$
converge. Therefore, if we replace $L$ by a finite extension of it,
$$ \rho_{1/d} \colon \Gal(\overline{L}/L) \ni \sigma \mapsto
   \exp \left( \frac{1}{d}\log(\rho(\sigma)) \right)
   \in \Q_p^{\times} $$
is a one dimensional $\Gal(\overline{L}/L)$-representation.
We denote it by $\Q_p(\frac{1}{d})$.
Then $\Q_p(\frac{1}{d})^{\otimes d} \cong \Q_p(1)$.
For $n \in \Z,\ n \geq 1$, we define
$\Q_p(\frac{n}{d}) := \Q_p(\frac{1}{d})^{\otimes n},
\ \Q_p(-\frac{n}{d}) := \Hom(\Q_p(\frac{n}{d}),\Q_p)$.
If $L$ is a finite extension of $\Q_p$,
we can similarly define $\C_p(a) \ (a \in \frac{1}{d}\Z)$
as in \S 4.3.

We define $K(l,S,\p)_{1/d}$ as follows.
Let notation be the same as in \S 5.2.
$K(l,S,\p)_{1/d}$ is an abelian group generated by
$\Gal(\overline{K}/K)$-representations $V$ satisfying
the following conditions modulo an equivalence relation $\sim$
as in Definition \ref{DefinitionGrothendieckGroups} :
\begin{enumerate}
\item The conditions 1, 2 in Definition
\ref{DefinitionGrothendieckGroups}
but we allow $k$ to be an element of $\frac{1}{d} \Z$
instead of $\Z$.
\item
Let $L$ be a finite extension of $K_{\p}$ such that
$\Q_p(\frac{1}{d})$ exists as
a $\Gal(\overline{L}/L)$-representation.
For $i,j \in \frac{1}{d}\Z$, we define
$$
h_{\p}^{i,\,j}(V) := \dim_L
(\Gr_{i+j}^W V \otimes_{\Q_p} \C_p(i))^{\Gal(\overline{L}/L)}.
$$
Then these numbers satisfy
$ \sum_{i, j \in \frac{1}{d}\Z} h_{\p}^{i,j}(V) = \dim_{\Q_p} V. $
It is easy to see that this condition is
independent of the choice of $L$.
\end{enumerate}

Similarly, we can define $K(l,S,\p,k)_{1/d}$ as in \S 5.2.
We have a direct sum decomposition as follows :
$$ K(l,S,\p)_{1/d} =
     \bigoplus_{k \in \frac{1}{d}\Z} K(l,S,\p,k)_{1/d}. $$
$K(l,S,\p)_{1/d}$ has a ring structure.
Moreover, we can define $p$-adic $E$-function
$E_{\p}(V;u,v) \in \Z[u^{1/d},v^{1/d}]$
for a $\Gal(\overline{K}/K)$-representation $V$
in $K(l,S,\p)_{1/d}$.
We can extend this to a ring homomorphism
$$ E_{\p} \colon K(l,S,\p)_{1/d} \longrightarrow
     \Z[u^{1/d},v^{1/d}]. $$

\begin{ex}
\label{ExampleFractionalTateTwists}
Assume that $\Q_p(\frac{1}{d})$ exists as
a $\Gal(\overline{K}/K)$-representation.
Note that this is satisfied if we replace $K$ by
a finite extension of it.
For $n \in \Z$,
$\Q_p(\frac{n}{d})$ is in $K(l,S,\p,-\frac{2n}{d})_{1/d}$
such that
$E_{\p}(\Q_p(\frac{n}{d});u,v) = u^{-n/d} v^{-n/d}$.
\end{ex}

\section{Conclusion --- the number of $\F_q$-rational points,
Galois representations, Hodge numbers}

We combine the results in \S 4 and \S 5.
Let notation be the same as in \S 5.2.

Let $X$ be a proper smooth variety over $K$
which has a proper smooth model $\X$ over $(\Spec \O_K) \backslash S$.
Then $H^{k}_{\text{\rm \'et}}(X_{\overline{K}},\Q_l)$ is
a $\Gal(\overline{K}/K)$-representation in $K(l,S,\p,k)$
by the Weil conjecture and $p$-adic Hodge theory (\S 4).
Let $X_{\overline{K}} = X \otimes_K \overline{K}$ and
$X_{\C} = X \otimes_K \C$.
We define a virtual representation
$$ [H^{\ast}_{\text{\rm \'et}}(X_{\overline{K}},\Q_l)]
  := \sum_{k} (-1)^k [H^{k}_{\text{\rm \'et}}(X_{\overline{K}},\Q_l)] $$
as an element in $K(l,S,\p)$.
Then we have the equality of two $E$-functions
$$ E(X_{\C};u,v)
     = E_{\p}([H^{\ast}_{\text{\rm \'et}}(X_{\overline{K}},\Q_l)];u,v) $$
by comparing the Hodge decomposition of $X_{\C}$
and the Hodge-Tate decomposition of $X_{K_{\p}} = X \otimes_K K_{\p}$.

\subsection{The proper smooth case}

By combining results in \S 4 and \S 5,
we have the following results
which connects the number of rational points and the Hodge numbers.

\begin{prop}
\label{RationalPointsGaloisRepresentations1}
Let $X$ be a proper smooth variety over $K$
which has a proper smooth model $\X$ over $(\Spec \O_K) \backslash S$.
Then we have
$$ |\X(\O_K/\p)| =
   \Tr(\Frob_{\p};[H^{\ast}_{\text{\rm \'et}}(X_{\overline{K}},\Q_l)])
      \qquad \text{for all} \quad \p \notin S. $$
\end{prop}

\begin{proof}
This follows from the Lefschetz trace formula for \'etale cohomology
as in \S 4.2.
\end{proof}

\begin{cor}
\label{RationalPointsGaloisRepresentationsCorollary1}
Let $X$ (resp. $Y$) be a proper smooth variety over $K$
which has a proper smooth model $\X$ (resp. $\Y$)
over $(\Spec \O_K) \backslash S$.
If $|\X(\O_K/\p)| = |\Y(\O_K/\p)|$ for all $\p \notin S$,
then
$[H^{\ast}_{\text{\rm \'et}}(X_{\overline{K}},\Q_l)]
= [H^{\ast}_{\text{\rm \'et}}(Y_{\overline{K}},\Q_l)]$
in $K(l,S,\p)$.
Therefore, we have
$$ E(X_{\C};u,v)
= E_{\p}([H^{\ast}_{\text{\rm \'et}}(X_{\overline{K}},\Q_l)];u,v)
= E_{\p}([H^{\ast}_{\text{\rm \'et}}(Y_{\overline{K}},\Q_l)];u,v)
= E(Y_{\C};u,v). $$
Namely, the Hodge numbers of $X_{\C}$ and $Y_{\C}$ are equal.
\end{cor}

\begin{proof}
By Proposition
\ref{RationalPointsGaloisRepresentations1},
we have
$$ \Tr(\Frob_{\p};[H^{\ast}_{\text{\rm \'et}}(X_{\overline{K}},\Q_l)])
  = \Tr(\Frob_{\p};[H^{\ast}_{\text{\rm \'et}}(Y_{\overline{K}},\Q_l)])
      \qquad \text{for all} \quad \p \notin S. $$
Hence we have
$[H^{\ast}_{\text{\rm \'et}}(X_{\overline{K}},\Q_l)]
 = [H^{\ast}_{\text{\rm \'et}}(Y_{\overline{K}},\Q_l)]$
in $K(l,S,\p)$ by Proposition \ref{CorChebotarev}.
\end{proof}

\subsection{A generalization --- the open smooth case}

Next we generalize
Proposition \ref{RationalPointsGaloisRepresentations1}
to open smooth varieties by a method of Deligne in
\cite{HodgeI},\cite{HodgeII}.

Let $X$ be a smooth variety over $K$ of dimension $n$
which is not necessarily proper.
Assume that there exists
a proper smooth variety $\overline{X} \supset X$ over $K$ such that
$\overline{X} \backslash X = \bigcup_{i=1}^{r} D_i$
is a SNCD on $\overline{X}$.
Let $I = \{ 1,\ldots,r \}$,\ $D_J = \bigcap_{j \in J} D_j$
for a nonempty subset $J \subset I$, and $D_{\emptyset} = \overline{X}$.
We consider the following formal sum
$$ [H^{\ast}_{c,\text{\rm \'et}}(X_{\overline{K}},\Q_l)]
 := \sum_{k}(-1)^k[H^{k}_{c,\text{\rm \'et}}(X_{\overline{K}},\Q_l)],$$
where $H^{k}_{c,\text{\rm \'et}}$ denotes \'etale cohomology
with compact support.
By the following Lemma \ref{LemmaOpenVarieties},
we see that the above is an equality in $K(l,S,\p)$.

\begin{lem}
\label{LemmaOpenVarieties}
Let $X$ be as above. Then we have
\begin{quote}
$\displaystyle E(X_{\C};u,v)
  = \sum_{J \subset I} (-1)^{|J|} \,E((D_J)_{\C};u,v)$
\quad $\in \Z[u,v]$, \\
$\displaystyle
[H^{\ast}_{c,\text{\rm \'et}}(X_{\overline{K}},\Q_l)]
  = \sum_{J \subset I} (-1)^{|J|}
       \,[H^{\ast}_{\text{\rm \'et}}((D_J)_{\overline{K}},\Q_l)]$
\quad $\in K(l,S,\p)$.
\end{quote}
Furthermore, we have the equality of two $E$-functions for $X$:
$$ E(X_{\C};u,v) =
   E_{\p}([H^{\ast}_{c,\text{\rm \'et}}(X_{\overline{K}},\Q_l)];u,v). $$
\end{lem}

\begin{proof}
Since $D_J$ is a proper smooth variety over $K$,
we have $E((D_J)_{\C};u,v) =
E_{\p}([H^{\ast}_{\text{\rm \'et}}((D_J)_{\overline{K}},\Q_l)];u,v)$.
Therefore, the second assertion immediately
follows from the first assertion.

We only prove the first assertion for $E(X_{\C};u,v)$
since we can prove the case of 
$[H^{\ast}_{c,\text{\rm \'et}}(X_{\overline{K}},\Q_l)]$
by the same way.
The Leray spectral sequence for
the inclusion $X_{\C} \hookrightarrow \overline{X}_{\C}$
induces a spectral sequence
$$ E_2^{i,j} = \bigoplus_{J \subset I\ \text{s.t.}\ |J| = j}
     H^i((D_J)_{\C},\Q)(-j)
     \Rightarrow H^{i+j}(X_{\C},\Q), $$
which defines the canonical mixed Hodge structure on
$H^k(X_{\C},\Q)$
(\cite{HodgeI},\cite{HodgeII}, 3.2).
For a finite dimensional $\Q$-vector space $V$ with
mixed Hodge structure, we define
$E(V;u,v) = \sum_{i,j} h^{i,j}(\Gr_{i+j}^W V)\,u^i v^j$.
Note that
$$ E(X_{\C};u,v) =
 \sum_{k=1}^{2n} (-1)^k E(H_c^k(X_{\C},\Q);u,v) $$
by definition.
By the above spectral sequence, we have
$$ \sum_{i,j} (-1)^{i+j} E(E_2^{i,j};u,v)
= \sum_{k=1}^{2n} (-1)^k E(H^k(X_{\C},\Q);u,v). $$
By Poincar\'e duality,
$H^k(X_{\C},\Q)$ is dual to $H_c^{2n - k}(X_{\C},\Q)(n)$.
Then we have
$$ E(H^k(X_{\C},\Q);u,v)
     = (uv)^n E(H_c^{2n-k}(X_{\C},\Q);u^{-1},v^{-1}). $$
On the other hand,
since $D_J$ is a proper smooth variety of dimension
$n - |J|$, $H^i((D_J)_{\C},\Q)(-|J|)$ is dual to
$H^{2n - 2|J| - i}((D_J)_{\C},\Q)(n)$ by Poincar\'e duality.
Hence we have
$$ E(E_2^{i,j};u,v) = \sum_{|J| = j}
   (uv)^n E(H^{2n - 2|J| - i}((D_J)_{\C},\Q);u^{-1},v^{-1}). $$
By combining them, we have Lemma \ref{LemmaOpenVarieties}.
\end{proof}

Next we consider the numbers of rational points valued
in finite fields.
For open smooth varieties,
we have the following generalization of
Proposition \ref{RationalPointsGaloisRepresentations1}.

\begin{prop}
\label{RationalPointsGaloisRepresentations2}
Let $X$ be a smooth variety over $K$.
Assume that there exist a proper smooth scheme
$\overline{\X}$ over $(\Spec \O_K) \backslash S$
and an open subscheme $\X \subset \overline{\X}$
whose generic fiber is $X$
such that
$\overline{\X} \backslash \X = \bigcup_{i=1}^{r} \D_i$
is a relative SNCD on
$\overline{\X} / (\Spec \O_K) \backslash S$
(see Definition \ref{DefinitionRelativeSNCD}).
Then we have
$$ |\X(\O_K/\p)| = \Tr(\Frob_{\p};
   [H^{\ast}_{c,\text{\rm \'et}}(X_{\overline{K}},\Q_l)])
      \qquad \text{for all} \quad \p \notin S. $$
\end{prop}

\begin{proof}
Let $I = \{ 1,\ldots,r \}$,
$\D_J = \bigcap_{j \in J} \D_j$
for a nonempty subset $J \subset I$,
and $\D_{\emptyset} = \overline{\X}$.
Then, by inclusion-exclusion principle, we have
$$
|\X(\O_K/\p)| = \sum_{J \subset I} (-1)^{|J|}\,|\D_J(\O_K/\p)|
\qquad \text{for all} \quad \p \notin S. $$
Since $\D_J$ is proper and smooth over $(\Spec \O_K) \backslash S$,
we have
$$
 |\D_J(\O_K/\p)| =
  \Tr(\Frob_{\p};[H^{\ast}_{\text{\rm \'et}}((D_J)_{\overline{K}},\Q_l)])
$$
by Proposition \ref{RationalPointsGaloisRepresentations1}.
On the other hand, we have
$$
 \Tr(\Frob_{\p};[H^{\ast}_{c,\text{\rm \'et}}(X_{\overline{K}},\Q_l)])
  = \sum_{J \subset I} (-1)^{|J|}
    \,\Tr(\Frob_{\p};[H^{\ast}_{\text{\rm \'et}}((D_J)_{\overline{K}},\Q_l)])
$$
by Lemma \ref{LemmaOpenVarieties}.
By combining them, we have Proposition
\ref{RationalPointsGaloisRepresentations2}.
\end{proof}

We note the following generalization of
Corollary \ref{RationalPointsGaloisRepresentationsCorollary1}
for open smooth varieties, although we don't use it later.

\begin{cor}
\label{RationalPointsGaloisRepresentationsCorollary2}
Let $X$ (resp. $Y$) be a smooth variety over $K$
satisfying the assumptions in
Proposition \ref{RationalPointsGaloisRepresentations2}.
Let $\X \subset \overline{\X}$ (resp. $\Y \subset \overline{\Y}$)
be a scheme over $(\Spec \O_K) \backslash S$
as in Proposition \ref{RationalPointsGaloisRepresentations2}.
If $|\X(\O_K/\p)| = |\Y(\O_K/\p)|$ for all $\p \notin S$,
then
$[H^{\ast}_{c,\text{\rm \'et}}(X_{\overline{K}},\Q_l)]
= [H^{\ast}_{c,\text{\rm \'et}}(Y_{\overline{K}},\Q_l)]$
in $K(l,S,\p)$.
Therefore, we have
$$ E(X_{\C};u,v)
= E_{\p}([H^{\ast}_{c,\text{\rm \'et}}(X_{\overline{K}},\Q_l)];u,v)
= E_{\p}([H^{\ast}_{c,\text{\rm \'et}}(Y_{\overline{K}},\Q_l)];u,v)
= E(Y_{\C};u,v). $$
\end{cor}

\begin{proof}
The proof is the same as
Corollary \ref{RationalPointsGaloisRepresentationsCorollary1}.
By Proposition
\ref{RationalPointsGaloisRepresentations2},
we have
$$ \Tr(\Frob_{\p};[H^{\ast}_{c,\text{\rm \'et}}(X_{\overline{K}},\Q_l)])
  = \Tr(\Frob_{\p};[H^{\ast}_{c,\text{\rm \'et}}(Y_{\overline{K}},\Q_l)])
      \qquad \text{for all} \quad \p \notin S. $$
Hence we have
$[H^{\ast}_{c,\text{\rm \'et}}(X_{\overline{K}},\Q_l)]
 = [H^{\ast}_{c,\text{\rm \'et}}(Y_{\overline{K}},\Q_l)]$
in $K(l,S,\p)$ by Proposition \ref{CorChebotarev}.
The equality of two $E$-functions
follows from Lemma \ref{LemmaOpenVarieties}.
\end{proof}

\section{Proof of the main theorem}

\begin{lem}
\label{LemmaSpecialization}
Let $f \colon \X \to T$ be a proper smooth morphism of schemes
of characteristic 0 and
$\D = \bigcup_{i=1}^{r} \D_i$ be
a relative SNCD on $\X/T$
(see Definition \ref{DefinitionRelativeSNCD}).
Assume that $T$ is connected.
Then, all fibers of $\X \backslash \D \to T$ have
the same $E$-functions defined in
Definition \ref{DefinitionE-function}.
\end{lem}

\begin{proof}
For a nonempty subset $J \subset \{ 1,\ldots,r \}$,
$\D_J = \bigcap_{j \in J} \D_j$ is proper and smooth over $T$.
Hence by the theorem of Deligne (\cite{Degen}, 5.5),
the Hodge numbers of all fibers of $\D_J \to T$ are the same.
On the other hand,
the $E$-function of a fiber of $\X \backslash \D \to T$
can be computed from the Hodge numbers of
a fiber of $\D_J \to T$ by Lemma \ref{LemmaOpenVarieties}.
Therefore, we have Lemma \ref{LemmaSpecialization}.
\end{proof}

\begin{proof}[Proof of Theorem \ref{MainTheorem}]
Let $\rho \colon Y \to X$ and $\rho' \colon Y' \to X$
be as in Theorem \ref{MainTheorem}.
Let $n$ be the dimension of $X$.
To avoid confusion, here
$E_{st}(X;u,v)_{\rho}$ (resp. $E_{st}(X;u,v)_{\rho'}$)
denotes the stringy $E$-function of $X$
defined by $\rho \colon Y \to X$ (resp. $\rho' \colon Y' \to X$)
as in Definition \ref{DefinitionStringyE-function}.
We shall prove the equality
$E_{st}(X;u,v)_{\rho} = E_{st}(X;u,v)_{\rho'}$.

\noindent
{\bf Step 1.}
\ Let $f \colon Y \dasharrow Y'$ be a birational map
between $Y$ and $Y'$ over $X$.
Let $Z$ be a resolution of singularities of
the closure of the graph of $f$ such that the exceptional
divisor of $Z \to X$ is a SNCD.
Let $\tau \colon Z \to Y$ and $\tau' \colon Z \to Y'$ be
a natural morphism.
Then we have the following commutative diagram.
$$
\xymatrix{
& Z \ar[dl]_{\tau} \ar[dr]^{\tau'} \\
Y \ar[dr]_{\rho} \ar@{-->}[rr]^{f} & & Y' \ar[dl]^{\rho'} \\
& X}
$$
If we prove $E_{st}(X;u,v)_{\rho} = E_{st}(X;u,v)_{\rho \circ \tau}$
and $E_{st}(X;u,v)_{\rho'} = E_{st}(X;u,v)_{\rho' \circ \tau'}$,
then we have Theorem \ref{MainTheorem}
since $\rho \circ \tau = \rho' \circ \tau'$.
In the followings, we only prove
$E_{st}(X;u,v)_{\rho} = E_{st}(X;u,v)_{\rho \circ \tau}$
since we can similarly prove
$E_{st}(X;u,v)_{\rho'} = E_{st}(X;u,v)_{\rho' \circ \tau'}$.

\noindent
{\bf Step 2.}
\ We shall show that we may assume everything is defined over
a number field.
Since $X,Y,Z,\rho,\tau$ are defined over a subfield $K'$ of $\C$
which is finitely generated over $\Q$.
Therefore, there exists an irreducible variety $T$ over
a number field $K$
such that the function field of $T$ is $K'$.
Furthermore, there exists a proper scheme $\widetilde{X}$,
proper smooth schemes $\widetilde{Y},\widetilde{Z}$,
and proper birational morphisms
$\tilde{\rho} \colon \widetilde{Y} \to \widetilde{X},
\ \tilde{\tau} \colon \widetilde{Z} \to \widetilde{Y}$
over $T$ such that
$\widetilde{X} \times_T \Spec \C = X,
\ \widetilde{Y} \times_T \Spec \C = Y,
\ \widetilde{Z} \times_T \Spec \C = Z,
\ \tilde{\rho} \times_T \Spec \C = \rho,
\ \tilde{\tau} \times_T \Spec \C = \tau$.
We write
$$ K_Y = \rho^{\ast} K_X + \sum_{i=1}^{r} a_i D_i,
\qquad
K_Z = (\rho \circ \tau)^{\ast} K_X + \sum_{j=1}^{s} b_j E_j $$
with $a_i \in \Q,\ a_i > -1,\ b_j \in \Q,\ b_j > -1$
(see Definition \ref{DefLogTerminal}).
Let $d$ be a positive integer such that
$(K_X)^{\otimes d}$ is a Cartier divisor on $X$.
Then $a_i,b_j$ are elements of $\frac{1}{d}\Z$.

By replacing $T$ by a Zariski open subset of it,
$(K_X)^{\otimes d}$ extends to a Cartier divisor
$(\Omega^n_{\widetilde{X}/T})^{\otimes d}$
on $\widetilde{X}$ over $T$, and we can write
$$ \Omega^n_{\widetilde{Y}/T}
    = \tilde{\rho}^{\ast} \Omega^n_{\widetilde{X}/T}
      + \sum_{i=1}^{r} a_i \widetilde{D}_i,
\qquad
\Omega^n_{\widetilde{Z}/T} = (\tilde{\rho} \circ \tilde{\tau})^{\ast}
\Omega^n_{\widetilde{X}/T} + \sum_{j=1}^{s} b_j \widetilde{E}_j, $$
where $\widetilde{D} = \bigcup_{i=1}^{r} \widetilde{D}_i$
(resp. $\widetilde{E} = \bigcup_{j=1}^{s} \widetilde{E}_j$) is
a relative SNCD on
$\widetilde{Y}/T$ (resp. $\widetilde{Z}/T$)
(see Definition \ref{DefinitionRelativeSNCD}).
By replacing $K$ by a finite extension of it,
there exists a $K$-rational point $t \in T(K)$.
Let $\widetilde{X}_t, \widetilde{Y}_t, \widetilde{Z}_t,
\tilde{\rho}_t, \tilde{\tau}_t$
be the fibers at $t$.
Then
$E_{st}(X;u,v)_{\rho} = E_{st}(\widetilde{X}_t;u,v)_{\tilde{\rho}_t}$
and $E_{st}(X;u,v)_{\rho \circ \tau}
= E_{st}(\widetilde{X}_t;u,v)_{\tilde{\rho}_t \circ \tilde{\tau}_t}$
by Lemma \ref{LemmaSpecialization}.
By replacing $X,Y,Z,\rho,\tau$ by
$\widetilde{X}_t, \widetilde{Y}_t, \widetilde{Z}_t,
\tilde{\rho}_t, \tilde{\tau}_t$,
we may assume $X,Y,Z,\rho,\tau$ are defined over a number field $K$.

\noindent
{\bf Step 3.}
\ By the same argument as above, there exists a finite set $S$ of
maximal ideals of $\O_K$, a proper scheme $\X$,
proper smooth schemes $\Y,\ZZ$,
and proper birational morphisms
$\bar{\rho} \colon \Y \to \X,\ \bar{\tau} \colon \ZZ \to \Y$
over $\T = (\Spec \O_K) \backslash S$
such that generic fibers of $\X,\Y,\ZZ,\bar{\rho},\bar{\tau}$
are $X,Y,Z,\rho,\tau$.
By enlarging $S$,
$(K_X)^{\otimes d}$ extends to a Cartier divisor
$(\Omega^n_{\X/\T})^{\otimes d}$ on $\X$ over $\T$
and we can write
$$ \Omega^n_{\Y/\T} = \bar{\rho}^{\ast} \Omega^n_{\X/\T}
      + \sum_{i=1}^{r} a_i \D_i,
\qquad
\Omega^n_{\ZZ/\T} = (\bar{\rho} \circ \bar{\tau})^{\ast}
\Omega^n_{\X/\T} + \sum_{j=1}^{s} b_j \E_j, $$
where $\D = \bigcup_{i=1}^{r} \D_i$
(resp. $\E = \bigcup_{j=1}^{s} \E_j$) is
a relative SNCD on $\Y/\T$ (resp. $\ZZ/\T$).

\noindent
{\bf Step 4.}
\ Here we compute $p$-adic integration.
Take a maximal ideal $\q \notin S$.
Let $K_{\q}$ be a $\q$-adic completion of $K$
and $\O_{K_{\q}}$ be the ring of integers of $K_{\q}$.
Let $q = |\O_{K}/\q|$ be the number of elements of
the residue field.

Let $\U_1,\ldots,\U_k$ be a finite open covering
of $\X$ over $\T$ such that
$(\Omega^n_{\X/\T})^{\otimes d}$ is a trivial line bundle
on each $\U_i$\ $(1 \leq i \leq k)$.
Let $\omega_i$ be a nowhere vanishing section of
$(\Omega^n_{\X/\T})^{\otimes d}$ on $\U_i$.

Then, by Proposition \ref{padicIntegration2}, we can compute
the $p$-adic integration of $\bar{\rho}^{\ast} \omega_i$
on $\bar{\rho}^{-1} \U_i(\O_{K_{\q}})$ as follows :
$$ \int_{\bar{\rho}^{-1} \U_i(\O_{K_{\q}})}
     |\bar{\rho}^{\ast} \omega_i|_p^{1/d}
     = \frac{1}{q^n} \sum_{J \subset \{ 1,\ldots,r \}}
         |(\D^{\circ}_J \cap \bar{\rho}^{-1} \U_i)(\F_q)|
         \prod_{j \in J} \frac{q-1}{q^{a_j+1}-1}, $$
where $\D^{\circ}_J$ are the same as in
Proposition \ref{padicIntegration2}.
Note that
$\text{div} (\bar{\rho}^{\ast} \omega_i) = \sum_{i=1}^{r} d a_i \D_i$
since $\omega_i$ is a nowhere vanishing section of
$(\Omega^n_{\X/\T})^{\otimes d}$.

Similarly, for $(\bar{\rho} \circ \bar{\tau})^{\ast} \omega_i$,
we have
$$ \int_{(\bar{\rho} \circ \bar{\tau})^{-1} \U_i(\O_{K_{\q}})}
     |(\bar{\rho} \circ \bar{\tau})^{\ast} \omega_i|_p^{1/d}
   = \frac{1}{q^n} \sum_{J' \subset \{ 1,\ldots,s \}}
      |(\E^{\circ}_{J'} \cap (\bar{\rho} \circ \bar{\tau})^{-1} \U_i)(\F_q)|
      \prod_{j' \in J'} \frac{q-1}{q^{b_{j'}+1}-1}, $$
where $\E^{\circ}_{J'}$ is the same as above.

On the other hand, by the change-of-variable formula for
$p$-adic integration, we have
$$ \int_{\bar{\rho}^{-1} \U_i(\O_{K_{\q}})}
   |\bar{\rho}^{\ast} \omega_i|_p^{1/d}
 = \int_{(\bar{\rho} \circ \bar{\tau})^{-1} \U_i(\O_{K_{\q}})}
   |(\bar{\rho} \circ \bar{\tau})^{\ast} \omega_i|_p^{1/d}. $$

Since the same is true for all finite intersections of $\U_i$,
by inclusion-exclusion principle, we conclude
$$ \frac{1}{q^n} \sum_{J \subset \{ 1,\ldots,r \}}
   |\D^{\circ}_J(\F_q)|
    \prod_{j \in J} \frac{q-1}{q^{a_j+1}-1}
     = \frac{1}{q^n} \sum_{J' \subset \{ 1,\ldots,s \}}
         |\E^{\circ}_{J'}(\F_q)|
         \prod_{j' \in J'} \frac{q-1}{q^{b_{j'}+1}-1}. $$

Note that the above argument works for every $\q \notin S$.

\noindent
{\bf Step 5.}
\ Fix a prime number $l = p$ and a maximal ideal $\p$ of $\O_K$
dividing $p$. 
By enlarging $S$, we may assume that $S$ contains
all maximal ideals of $\O_K$ dividing $p$.
We shall work on the level of the Grothendieck group
$K(l,S,\p)_{1/d}$ of $\Gal(\overline{K}/K)$-representations
introduced in \S 5.

We rewrite the conclusion of Step 4 in the following form
by multipling
$q^n \cdot \prod_{j=1}^{r} (q^{a_j+1}-1) \cdot
\prod_{j'=1}^{s} (q^{b_{j'}+1}-1)$ on both sides :
\begin{align*}
\prod_{j'=1}^{s} (q^{b_{j'}+1}-1)
   \sum_{J \subset \{ 1,\ldots,r \}} \bigg(
   |\D^{\circ}_J(\F_q)|
   \prod_{j \in J} (q-1) \prod_{j \notin J} (q^{a_j+1}-1)
   \bigg) \\
&\hspace*{-3.7in}= \prod_{j=1}^{r} (q^{a_j+1}-1)
   \sum_{J' \subset \{ 1,\ldots,s \}} \bigg(
   |\E^{\circ}_{J'}(\F_q)|
   \prod_{j' \in J'} (q-1) \prod_{j' \notin J'} (q^{b_{j'}+1}-1)
   \bigg).
\end{align*}

By replacing $K$ by a finite extension of it, we may assume
that $\Q_l(\frac{1}{d})$ exists as
a $\Gal(\overline{K}/K)$-representation
(see Example \ref{ExampleFractionalTateTwists}).
Recall that the image of $\Frob_{\q}$ in $\Gal(\overline{\F}_q/\F_q)$
is the inverse of the $q$-th power automorphism $x \mapsto x^q$
of $\overline{\F}_q$.
Therefore, for $m \in \frac{1}{d}\Z$,
$\Tr(\Frob_{\q};\Q_l(m)) = q^{-m}$.

Hence, we have the following equality in $K(l,S,\p)_{1/d}$ :
\begin{quote}
\hspace*{-0.3in} $\displaystyle
   \prod_{j'=1}^{s} ([\Q_l(-b_{j'}-1)]-1)
   \sum_{J \subset \{ 1,\ldots,r \}}
   \bigg(
   [H^{\ast}_{c,\text{\rm \'et}}((D^{\circ}_J)_{\overline{K}},\Q_l)]
   \prod_{j \in J} ([\Q_l(-1)]-1)$ \\
\hspace*{0.6in} $\displaystyle
   \cdot 
   \prod_{j \notin J} ([\Q_l(-a_j-1)]-1)
   \bigg)$ \\
\hspace*{-0.3in} $\displaystyle
  = \prod_{j=1}^{r} ([\Q_l(-a_j-1)]-1)
   \sum_{J' \subset \{ 1,\ldots,s \}}
   \bigg(
   [H^{\ast}_{c,\text{\rm \'et}}((E^{\circ}_{J'})_{\overline{K}},\Q_l)]
   \prod_{j' \in J'} ([\Q_l(-1)]-1)$ \\
\hspace*{0.8in} $\displaystyle
   \cdot \prod_{j' \notin J'} ([\Q_l(-b_{j'}-1)]-1)
   \bigg)$,
\end{quote}
since the traces of $\Frob_{\q}$ on the both sides
are equal for all $\q \notin S$
(see Proposition \ref{RationalPointsGaloisRepresentations2}
and Proposition \ref{CorChebotarev}).
Note that $1 \in K(l,S,\p)_{1/d}$ denotes the class of
the trivial $\Gal(\overline{K}/K)$-representation.

Since $E_{\p}(\Q_l(m);u,v) = u^{-m} v^{-m}$ for $m \in \frac{1}{d}\Z$
by Example \ref{ExampleFractionalTateTwists},
we have
\begin{align*}
   \prod_{j'=1}^{s}((uv)^{b_{j'}+1}-1)
   \sum_{J \subset \{ 1,\ldots,r \}}
   \bigg(
   E((D^{\circ}_J)_{\C};u,v)
   \prod_{j \in J} (uv-1) \prod_{j \notin J} ((uv)^{a_j+1}-1)
   \bigg) \\
&\hspace*{-5in}=
   \prod_{j=1}^{r}((uv)^{a_j+1}-1)
   \sum_{J' \subset \{ 1,\ldots,s \}}
   \bigg(
   E((E^{\circ}_{J'})_{\C};u,v)
   \prod_{j' \in J'} (uv-1) \prod_{j' \notin J'} ((uv)^{b_{j'}+1}-1)
   \bigg).
\end{align*}
by Lemma \ref{LemmaOpenVarieties}.
By Definition \ref{DefinitionStringyE-function},
this proves
$E_{st}(X_{\C};u,v)_{\rho} = E_{st}(X_{\C};u,v)_{\rho \circ \tau}$
and hence Theorem \ref{MainTheorem}.
\end{proof}

\begin{rem}
If we take an appropriate $\p$ in Step 5,
we can use the result of Fontaine-Messing
(Remark \ref{Remarkp-adicHodgeTheory}, \cite{FontaineMessing}).
Therefore, we don't need the full version of
the Hodge-Tate decomposition for Theorem \ref{MainTheorem}.
\end{rem}

\end{document}